\documentstyle[12pt]{article}

\begin{document}

\noindent Title: On the uniphase solutions of 
the nonlinear damped wave equation


\noindent Authors: {\bf S. Bir\u{a}ua\c{s}, D. Opri\c{s}} Departament of Mathematics 
The West University of Timi\c{s}oara Timi\c{s}oara 1900, Romania\\
\noindent Comments: 19 pages, Tex\\
\noindent Subj-class: Partial Differential Equations\\
\noindent MSC-class: 39A13 (Primary) 35B10 (Secondary)\\
In this paper we study the steady uniphase and multiphase 
solutiuons of the discretized nonlinear damped wave equation. Conditions for the
stability and instability of the steady solutions are given; in the instability
case the linear stable and unstable manifolds are described.
\vspace*{1cm}

\begin{center}
\Large \bf 1. Introduction
\end{center}

\vspace*{.5cm}
We consider the following mixed problem:
\begin{equation}
u_{tt}=(\sigma(u_{x}))_{x}+\varepsilon u_{xxt}, \; x\in (0,1),\; t>0,\;
\varepsilon \in {\bf R}
\end{equation}
\begin{equation}
u(0,t)=0,\; u(1,t)=P,\; t>0
\end{equation}
\begin{equation}
u(x,0)=u_{0}(x),\; u_{t}(x,0)=u_{1}(x),\; x\in (0,1)
\end{equation}
where $P>0$ and $\sigma :{\bf R} \to {\bf R}$ is an application with the following
properties: $\sigma (0)=0$, $\sigma$ is continously differentiable on ${\bf R}$,
$\sigma (\xi) >0$, $\forall \xi >0$, $\sigma (\xi )<0$, $\forall \xi <0$; there
exist $\overline{\alpha}$, $\underline{\beta}$ such that $0<\overline{\alpha}<
\underline{\beta}$ and $\sigma^{\prime}(\xi )>0$, $\forall \xi \in (-\infty ,
\overline{\alpha}) \cup (\underline{\beta} ,\infty )$, $\sigma^{\prime}(\xi)<0$,
$\forall \xi \in (\overline{\alpha}, \underline{\beta})$, and $\underline{\alpha}, \overline{\beta},$
such that $0<\underline{\alpha}<\overline{\beta}<\underline{\beta}<\overline{\beta},$
such that $\sigma(\underline{\alpha})=\sigma(\underline{\beta})$ and
$\sigma(\overline{\alpha})=\sigma(\overline{\beta}).$

In the paper [1] it is shown that for $\varepsilon=1$ the mixed problem 1-3 has a
unique solution on the appropiate function space. The uniphase and multiphase steady
solutions are defined and it is studied the stability of these steady solutions.
It is presented also a discretization in the spatial variable $x$ of the equation
({\bf 1}.1) and some numerical simulations are presented. For $\varepsilon =0$ the problem
was studied in [2].

For $\sigma (\xi )=a^{2}\xi$, $a\in {\bf R}$ we obtain the wave equation studied
in a lot of papers.

The aim of this paper is the study of uniphase steady solutions and partially
the study of multiphase steady solutions of the system obtained via a discretization
of ({\bf 1}.1) in spatial variable $x$ and then a discretization in temporal variable $t$.

In the section {\bf 2} it is shown that the system obtained from the discretization of
({\bf 1}.1) in the spatial variable $x$ represents an equation Euler-Lagrange by rapport
with a Lagrange function discrete continuous with damped term.

The system has a finite number of steady solutions (uniphase and multiphase).

In the section {\bf 3}, using the linearized system of the system used in the section
{\bf 2}, we study the stability of uniphase steady solutions. We proove that such
a solution is stable if $P\in (\underline{\alpha}, \overline{\alpha})$ or $P\in
(\underline{\beta}, \overline{\beta})$ and $\varepsilon >0$. If $P\in (\underline{\alpha},
\underline{\beta})$ the uniphase steady solution is hiperbolic so it is unstable.
We have determined a polycycle curve for this solution.

For $P\in (\underline{\alpha},\overline{\alpha})$ or $P\in (\underline{\beta},
\overline{\beta})$, $\varepsilon =0$, we find a curve $u^{c}(t)$, $t\in {\bf R}$
such that every of this components is periodic with the period
$$T_{k}=\frac{\pi}{n\sqrt{\tau}} cosec\; \frac{k\pi}{2n}$$

For the steady solution 2-phase we present stability conditions depending on $\varepsilon$,
$\tau_{1}=\sigma^{\prime}(\alpha)$, $\rho=\sigma^{\prime}(\beta)$, $\alpha \in
(\underline{\alpha},\overline{\alpha})$, $\beta \in (\underline{\beta},\overline{\beta})$.

In the section {\bf 4} it is shown that the discretization of the equation ({\bf 1}.1)
leads us to a difference system of equations which represents the equations Euler-Lagrange
with damping for an associated discrete Lagrange function. The steady solutions,
uniphase and multiphase, are presented and we proove that they are in finite
number.

Also, in the section {\bf 4} we are studing the stability of the uniphase solutions,
the conditions for a uniphase steady solution to be hiperbolic and the linear stable
and unstable associated manifolds.

A similar study for the multiphase steady solutions will be done in a future paper.

\setcounter{equation}{0}
\begin{center}
\Large \bf 2. Semidiscretization of the equation $u_{tt}=(\sigma (u_{x}))_{x}+\varepsilon
u_{xxt}$
\end{center}

\vspace*{.5cm}

For $\varepsilon =0$, the equation ({\bf 1}.1) represents the equation Euler-Lagrange
for the Lagrange function $L:J^{2}({\bf R}^{2}, {\bf R}) \to {\bf R}$,
$L(x,t,u,u_{x},u_{tx})=\frac{1}{2} u^{2}_{x}-\sigma (u_{x})$, where $J^{2}({\bf R}^{2},
{\bf R})$ is the bundle of the jets of order 2 of the fibrate $\pi :{\bf R}^{3}
\to{\bf R}^{2}$.

Be $x_{k}=kh_{1}$, $k=1,2,...,n$ the division points of the interval $[0,1]$,
$h_{1}=1/n$ and $u_{k}(t)=u(kh_{1},t)$, $k=1,2,...,n$. For the boundary conditions
we consider: $u_{0}=u(0,t)=0$, $u_{n}=u(1,t)=P$.

The partial derivatives $u_{x}$, $u_{xx}$, $(\sigma (u_{x}))_{x}$ will be approximate by:
$$u_{x}(kh_{1},t)\sim \frac{1}{h_{1}} (u_{k}(t)-u_{k-1}(t))$$
\begin{equation}
u_{xx}(kh_{1},t)\sim \frac{1}{h_{1}^{2}} (u_{k+1}(t)-2u_{k}(t)+u_{k-1}(t))
\end{equation}
$$(\sigma(u_{x}))_{x}(kh_{1},t)\sim \frac{1}{h_{1}} \left[ \sigma \left( \frac{1}{h_{1}}
(u_{k+1}(t)-u_{k}(t)\right) -\sigma \left( \frac{1}{h_{1}} (u_{k}(t)-u_{k-1}(t)
\right) \right] $$

Let be the sequence $(u_{k},\dot{u}_{k})\in T{\bf R}$, $k=1,...,n$ on the tangent bundle
at $\bf R$ and $L:T{\bf R}\to {\bf R}$ the Lagrange function defined by:
\begin{equation}
L(u_{k},\dot{u}_{k})=\frac{1}{2}\dot{u}^{2}_{k}+w \left( \frac{1}{h_{1}}
(u_{k}-u_{k-1})\right) \; ,\; k=1,...,n-1
\end{equation}
where $w :{\bf R}\to {\bf R}$ with $w^{\prime} (\xi)=\sigma (\xi)$. The
action of $L$ is defined by:
\begin{equation}
{\cal{A}}(u,\dot{u})=\sum_{k=1}^{n} L(u_{k},\dot{u}_{k})
\end{equation}

In order to obtain the first variation of $\cal{A}$ we consider the sequence
$(u_{k}(\eta),\dot{u}_{k}(\eta))\in T{\bf R}$, $k=1,...,n$ with $\eta \in (-a,a)$ and
$u_{k}(0)=u_{k}$, $\dot{u}_{k}(0)=\dot{u}_{k}$. The action ({\bf 2}.3) on this
sequence is:
$${\cal{A}} (\eta)=\sum_{k=1}^{n-1} L(u_{k}(\eta),\dot{u}_{k}(\eta)) $$

The first variation of ({\bf 2}.3) is given by
$$\left. \frac{\partial {\cal{A}} (\eta)}{\partial \eta} \right|_{\eta=0} =0 $$

The first variation of ({\bf 2}.3) is:
\begin{equation}
\ddot{u}_{k}-\frac{1}{h_{1}} \sigma (\Delta u_{k})+\frac{1}{h_{1}} \sigma (\Delta
u_{k-1})=0\; ,\; k=1,...,n-1
\end{equation}
where
$$\Delta u_{k+1}=\frac{1}{h_{1}} (u_{k+1}-u_{k})\; ,\; k=1,...,n-1$$
and
$$\Delta u_{n}=\frac{1}{h_{1}} (P-u_{n-1})$$

The system ({\bf 2}.4) represent the semidiscretized system of the equation ({\bf 1}.1)
for $\varepsilon =0$.

The semidiscretized system of the equation ({\bf 1}.1) for $\varepsilon \neq 0$ is:
\begin{equation}
\ddot{u}_{k}-\frac{1}{h_{1}} \sigma (\Delta u_{k+1})+\frac{1}{h_{1}} \sigma
(\Delta u_{k})=\frac{\varepsilon}{h_{1}^{2}} (\dot{u}_{k+1}-2\dot{u}_{k}+\dot{u}_{k-1})\;
, \; k=1,...,n-1
\end{equation}
and represents the equation Euler-Lagrange for the function $L$ defined in ({\bf 2}.2)
with the disperssion:
$$\frac{\varepsilon}{h_{1}^{2}} (\dot{u}_{k+1}-2\dot{u}_{k}+\dot{u}_{k-1}) $$

To the system ({\bf 2}.5) we associate the equivalent system on $T^{*}R$:
$$\dot{u}_{k}=v_{k}\; ,\; \; k=1,...,n-1 $$
\begin{equation}
\dot{v}_{k}=\frac{1}{h_{1}} \left[ \sigma (\Delta u_{k+1})-\sigma (\Delta u_{k}) 
\right] +\frac{\varepsilon}{h_{1}^{2}} (v_{k+1}-2v_{k}+v_{k-1})
\end{equation}
with the action:
\begin{equation}
V(u,v)=\sum_{k=1}^{n-1} \left[ \frac{1}{2}v^{2}_{k}+w(\Delta u_{k})\right] +w
(\Delta u_{n})
\end{equation}

From ({\bf 2}.3) and ({\bf 2}.5) it follows that:
\begin{equation}
\frac{d{\cal{A}}(u,\dot{u}}{dt}) =-\sum_{k=1}^{n} \left[ \frac{1}{h_{1}} (\dot{u}_{k+1}
-\dot{u}_{k}) \right]^{2}
\end{equation}

From ({\bf 2}.5) and ({\bf 2}.6) it follows that any steady solution of the system
({\bf 2}.5) satisfies the conditions:
\begin{equation}
\sigma(\Delta u_{k})=C,\; k=1,...,n,\; \; \; \sum_{k=1}^{n} \frac{1}{h_{1}} \Delta u_{k}
=nP, \; C\in {\bf R}
\end{equation}

Because the function $\sigma$ is not monotone we obtain two types of steady solutions:

\underline{\bf Definition 2.1.} A uniphase steady solution $\bar{u}$ for the
system ({\bf 2}.5) is a steady solution with the property: $\bar{u}_{k}=kh_{1}P$,
$k=1,...,n-1$.

\underline{\bf Definition 2.2.} A multiphase steady solution for the system ({\bf 2}.5)
is a steady solution with the property: $u_{k+1}-u_{k}=h_{1}f(k)$, $k=1,...,n-1$.

We denote by $E$ the set of multiphase steady solutions and:
$$E^{+}=\{ \bar{u}\in E : \sigma^{\prime} (\Delta u_{k})>0,\; k=1,...,n \} $$
$$E^{-}=\{ \bar{u}\in E : \exists k\in \{ 1,...,n\} :\sigma^{\prime} (\Delta u_{k})<0 \} $$

In [1] it is shown that:

\underline{\bf Lemma 2.1.} For any steady solution $\bar{u}\in E^{+}$ there is
$\alpha \in (\underline{\alpha},\bar{\alpha})$ and $\beta \in (\underline{\beta},
\bar{\beta})$, $C\in (\underline{\sigma}, \bar{\sigma})$, such that:

1) $\sigma (\Delta \bar{u}_{k})=C,\; \; k=1,...,n$

2) $\Delta (\bar{u}_{k})=\alpha$ or $\Delta (\bar{u}_{k})=\beta , \; k=1,...,n$

3) $k\alpha +(n-k)\beta =nP$

\underline{\bf Lemma 2.2.} If $P<\underline{\alpha}$ or $P>\bar{\beta}$ then
$E^{+}=\phi$ 

\underline{\bf Lemma 2.3.} The system ({\bf 2}.5) has a finite number of steady solutions.

\newpage
\setcounter{equation}{0}
\begin{center}
\Large \bf 3. Properties of the steady solutions of the system (2.5)
\end{center}

\vspace*{.5cm}

The general form of the system ({\bf 2}.6) is:
\begin{equation}
F_{k}(u_{k},u_{k+1},\dot{u}_{k-1},\dot{u}_{k},\dot{u}_{k+1},\ddot{u}_{k})=0,\; 
k=1,2,...,n
\end{equation}

Let $\bar{u}$ be a steady solution of the system ({\bf 2}.6). The linearized system
in a neighborhood of $\bar{u}$ is given by:
\begin{equation}
\left. 
\frac{\partial F_{k} \left( \bar{u} +\displaystyle \sum_{j=1}^{n} \eta^{j}w_{j}, \dot{\bar{u}}+\sum_{j=1}^{n}
\eta^{j}\dot{w}_{j}, \ddot{\bar{u}}+\sum_{j=1}^{n}\eta^{j}\ddot{w}_{j} \right) }{\partial
\eta^{j}} \right|_{\eta^{j}} =0
\end{equation}
$\hspace*{10cm} j=1,...,n$ \\
where $\eta_{j}\in (-a,a),\; j=1,...,n$.

The linearized system associated to ({\bf 2}.5) is:
\begin{equation}
\ddot{w}_{j}-\frac{1}{h_{1}^{2}} \sigma^{\prime} (\Delta \bar{u}_{j+1})(w_{j+1}-w_{j})+
\frac{1}{h_{1}^{2}} \sigma^{\prime} (\Delta \bar{u}_{j})(w_{j}-w_{j-1})=
\end{equation}
$$=\frac{\varepsilon}{h_{1}^{2}} (\dot{w}_{j+1}-2\dot{w}_{j}+\dot{w}_{j-1}),\; \;
j=1,...,n-1$$

A solution for the system ({\bf 3}.3) is:
\begin{equation}
w_{j}=exp\; (\lambda j)exp\; (ia_{k}j),\; a_{k}=\frac{\pi k}{n},\; k=1,...,n-1,
\; \lambda \in {\bf C}
\end{equation}

From ({\bf 3}.3) and ({\bf 3}.5) yields:
$$\lambda^{2}exp\; ia_{k}-\frac{1}{h_{1}^{2}} \sigma^{\prime} (\Delta \bar{u}_{k+1})
exp\; ia_{k}(exp\; ia_{k-1})+\frac{1}{h_{1}^{2}} \sigma^{\prime} (\Delta \bar{u}_{k})
(exp\; ia_{k-1}) =$$
\begin{equation}
=\frac{\varepsilon}{h_{1}^{2}} \lambda (exp\; ia_{k-1})^{2}
\end{equation}

\underline{\bf Theorem 3.1.} If $P\in (\underline{\alpha},\bar{\alpha})$ or $P\in
(\underline{\beta},\bar{\beta})$ the uniphase steady solution $\bar{u}_{k}=kh_{1}P$,
$k=1,...,n-1$ of the system ({\bf 2}.5) is asymptotically stable if $\varepsilon >0$
and unstable if $\varepsilon <0$.

\underline{\it Proof: } For the uniphase steady solution $\bar{u}$, the equations
({\bf 3}.5) become:
\begin{equation}
\lambda^{2} exp\; ia_{k}-\frac{1}{h_{1}^{2}} \tau (exp\; ia_{k-1})^{2}-\frac{\varepsilon}
{h_{1}^{2}} \lambda (exp\; ia_{k-1})^{2}=0,\; k=1,...,n-1
\end{equation}
where $\tau =\sigma^{\prime}(P)$. We denote:
$$\mu_{k}=\frac{(exp\; ia_{k-1})^{2}}{exp\; ia_{k}}=-4sin^{2}\; \frac{\pi k}{2n}$$

The equations ({\bf 3}.6) become:
\begin{equation}
\lambda^{2}-\varepsilon \mu_{k}n^{2}\lambda -\mu_{k}n^{2}\tau=0,\; k=1,...,n-1
\end{equation}
with the roots:
\begin{equation}
\lambda^{2}=\frac{1}{2} \mu_{k}n^{2} \left[ \varepsilon \pm \left( \varepsilon^{2}
+\frac{4\tau}{\mu_{k}n^{2}}\right)^{\frac{1}{2}} \right],\; k=1,...,n-1
\end{equation}

If $\varepsilon >0$, because $\mu_{k}<0$ yields that $Re\; \lambda_{k} <0$ so it
follows that the steady uniphase solution $\bar{u}$ is asymptotically stable.

If $\varepsilon <0$, it follows that $Re\; \lambda_{k} <0$ so the steady uniphase
solution is unstable.

\underline{\bf Theorem 3.2.} If $P\in (\bar{\alpha}, \underline{\beta})$ the uniphase
steady solution $\bar{u}$ of the system ({\bf 2}.5) is hyperbolic so it is unstable.

\underline{\it Proof: } For the uniphase steady solution $\bar{u}$ in the equations
({\bf 3}.7) become:
\begin{equation}
\lambda^{2}-\varepsilon \mu_{k}n^{2}\lambda +\mu_{k}n^{2}\rho =0,\; k=1,...,n-1
\end{equation}
where $\rho=-\sigma^{\prime}(P)$; the solutions of these equations are:
\begin{equation}
\lambda^{+}_{k}=\frac{1}{2}\mu_{k}n^{2}\left[ \varepsilon +\left( \varepsilon^{2}
-\frac{4\rho}{\mu_{k}n^{2}} \right)^{\frac{1}{2}} \right] 
\end{equation}
$$\lambda^{-}_{k}=\frac{1}{2}\mu_{k}n^{2}\left[ \varepsilon -\left( \varepsilon^{2}
-\frac{4\rho}{\mu_{k}n^{2}} \right)^{\frac{1}{2}} \right] ,\; k=1,...,n-1 $$

The real solutions $\lambda^{+}$, $\lambda^{-}$ have the property: $\lambda^{+}
\cdot \lambda^{-}=\mu_{k}n^{2}\rho <0$. It follows that they have opposite signs.
The steady solution is unstable, it is hyperbolic.

The linear manifolds associated to the steady solution of hyperbolic type are:
\begin{equation}
\tilde{u}_{k}^{u}(t)=h_{1}kP+exp\; (t\lambda_{k}^{+}) \sum_{j=1}^{n-1} \eta_{j}
exp\; (ia_{k}j),\; t\in{\bf R}, \; k=1,...,n-1
\end{equation}
\begin{equation}
\tilde{u}_{k}^{s}(t)=h_{1}kP+exp\; (t\lambda_{k}^{-}) \sum_{j=1}^{n-1} \eta_{j}
exp\; (ia_{k}j),\; t\in{\bf R}, \; k=1,...,n-1
\end{equation}

$\tilde{u}_{k}^{u}$ and $\tilde{u}_{k}^{s}$ have the properties:
$$\lim_{t\to -\infty} \tilde{u}_{k}^{u}(t)=\bar{u}_{k},\; \lim_{t\to \infty}
\tilde{u}_{k}^{s}(t)=\bar{u}_{k},\; k=1,...,n-1 $$

\underline{\bf Theorem 3.3.} The curves $u_{kk+1}^{h}$, $k=1,2,...,n$, $t\in {\bf R}$
given by:
\begin{equation}
u_{kk+1}^{h}=\frac{\tilde{u}_{k}^{s} (t)}{\tilde{u}_{k+1}^{s} (t)}\bar{u}_{k+1}+
\frac{\tilde{u}_{k+1}^{u} (t)}{\tilde{u}_{k}^{u} (t)}\bar{u}_{k},\; t\in {\bf R}
\end{equation}
are policycles for the uniphase steady solution, that is
$$\lim_{t\to \infty} u^{h}_{kk+1} (t)=\bar{u}_{k}$$

\underline{\it Proof: } From ({\bf 3}.12) and ({\bf 3}.13) we have:
$$\lim_{t\to \infty} u^{h}_{kk+1} (t)=\frac{\bar{u}_{k}}{\bar{u}_{k+1}}
\bar{u}_{k+1}+\bar{u}_{k}\lim_{t\to \infty} \frac{\tilde{u}_{k+1}^{u}(t)}
{\tilde{u}_{k}^{u}(t)}$$

From ({\bf 3}.11) yields:
$$\lim_{t\to \infty} \frac{\tilde{u}_{k+1}^{u}(t)}{\tilde{u}_{k}^{u}(t)}=
\lim_{t\to \infty} \frac{h_{1}(k+1)P+\displaystyle exp(t\lambda_{k+1}^{+})\sum_{j=1}^{n}
\eta_{j}exp(ia_{k+1}j)}{h_{1}(k)P+\displaystyle exp(t\lambda_{k}^{+})\sum_{j=1}^{n} \eta_{j}exp(ia_{k}j)}=$$
$$=\frac{\lambda^{+}_{k+1}}{\lambda^{+}_{k}} \lim_{t\to \infty} exp(t(\lambda_{k+1}^{+}-
\lambda^{+}_{k})) \frac{\displaystyle \sum_{j=1}^{n} \eta_{j}exp(ia_{k+1}j)}
{\displaystyle \sum_{j=1}^{n} \eta_{j}exp(ia_{k}j)}$$

From ({\bf 3}.10) we have $\lambda^{+}_{k+1}<\lambda^{+}_{k}$ and $\lambda^{-}_{k+1}
<\lambda^{-}_{k}$. It follows that
$$\lim_{t\to \infty} \frac{\tilde{u}_{k+1}^{h}(t)}{\tilde{u}_{k}^{h}(t)}=0$$

Thus we obtain:
$$\lim_{t\to \infty} u^{h}_{k+1} (t)=\bar{u}_{k}$$

Using ({\bf 3}.11) and ({\bf 3}.13) it follows:
$$\lim_{t\to -\infty} u^{h}_{kk+1} (t)=\bar{u}_{k+1}+\bar{u}_{k+1}
\lim_{t\to -\infty} \frac{\tilde{u}_{k}^{s}(t)}{\tilde{u}_{k+1}^{s}(t)}$$

From ({\bf 3}.12) yields:
$$\lim_{t\to -\infty} \frac{\tilde{u}_{k}^{s}(t)}{\tilde{u}_{k+1}^{s}(t)}=
\frac{\lambda^{-}_{k}}{\lambda^{-}_{k+1}} \lim_{t\to -\infty} exp\; (t(\lambda^{-}_{k}
-\lambda^{-}_{k+1})\cdot $$
$$\cdot \frac{\displaystyle \sum_{j=1}^{n} \eta_{j}exp(ia_{k}j)}{\displaystyle 
\sum_{j=1}^{n} \eta_{j}exp(ia_{k+1}j)}=0$$

Thus we obtain:
$$\lim_{t\to -\infty} u_{kk+1}^{h}=\bar{u}_{k+1}$$

The curves ({\bf 3}.13) are a polycycle considering:
$$u_{n+1}=u_{1}$$

\underline{\bf Theorem 3.4.} If $P\in(\underline{\alpha},\bar{\alpha})$ or
$P\in(\underline{\beta},\bar{\beta})$ and $\varepsilon =0$ the curves:
\begin{equation}
u_{k}^{c}(t)=\bar{u}_{k}+exp\; (\lambda^{c}t)\sum_{j=1}^{n} \eta_{j}exp\; (ia_{k}j),
\; k=1,...,n,\; t\in {\bf R}
\end{equation}
are periodic, with the period $\displaystyle T_{k}=\frac{\pi}{n\sqrt{\tau}} sec\;
\frac{k\pi}{2n}$, $k=1,...,n$, where $\lambda^{c}=\pm 2in\sqrt{\tau} sin\; \displaystyle
\frac{k\pi}{2n}$.

\underline{\it Proof: } For $\varepsilon =0$ the equation ({\bf 3}.7) is:
$$\lambda^{2}-\mu_{k}n^{2}\tau =0,\; k=1,...,n-1$$
with the roots:
$$\lambda^{c}=\pm 2in\sqrt{\tau} sin \; \frac{k\pi}{2n}$$

The uniphase steady solution $\bar{u}_{k}=kh_{1}P$ is center for the system ({\bf 2}.5)
with $\varepsilon=0$. From the condition $u_{k}^{c}(t)=u_{k+1}^{c}(t+T_{k})$,
$k=1,...,n-1$, it follows that $exp\; (\lambda^{c}(t+T_{k}))=exp\; (\lambda^{c}(t)$.
Hence: $2n\sqrt{\tau} sin \displaystyle \frac{k\pi}{2n} T_{k}=2\pi$ or
$$T_{k}=\frac{\pi}{n\sqrt{\tau}} cosec \frac{k\pi}{2n}$$

We continue our study with a multiphase steady solution $\bar{u}\in E^{+}$. To
be more specific we will consider only 2-phase steady solutions. From {\it Lemma
2.1.} it follows that for $\alpha \in (\underline{\alpha},\bar{\alpha})$, $\beta
\in (\underline{\beta},\bar{\beta})$, with $I(\alpha)=\{ i\in \{ 1,2,...,n\} $,
$\Delta \bar{u}_{i}=\alpha \} $, $I(\beta)=\{ i\in \{ 1,...,n\} ,\Delta \bar{u}_{i}
=\beta \} $ we have
$$I(\alpha) \cup I(\beta)=\{ 1,...,n\}$$
$$card\; I(\alpha)\cdot \alpha +(n-card\; I(\beta))\beta =nP$$

Let consider the following sets:
$$J(\alpha)=\{ i\in I(\alpha):\; i+1\in I(\alpha) \} $$
$$J(\beta)=\{ i\in I(\beta):\; i+1\in I(\beta) \} $$
and
$$K(\alpha ,\beta)=\{ i\in I(\alpha), \; i+1\in I(\beta) \} $$

From ({\bf 3}.5) we obtain:
\begin{equation}
\lambda^{2}+4\varepsilon n^{2}sin^{2}\frac{a_{p}}{2}\lambda +4n^{2}sin^{2}\frac{a_{p}}
{2}\tau =0, \; p\in J(\alpha)
\end{equation}
\begin{equation}
\lambda^{2}+4\varepsilon n^{2}\lambda sin^{2}\frac{a_{p}}{2}+4\rho n^{2}sin^{2}
\frac{a_{p}}{2}\tau =0, \; p\in J(\beta)
\end{equation}
\begin{equation}
\lambda^{2}+4\varepsilon n^{2}\lambda sin^{2}\frac{a_{p}}{2}-n^{2}[\rho (exp\; ia_{p-1})+
\tau (exp\; (-ia_{p})-1 ]=0,\; p\in K(\alpha ,\beta)
\end{equation}
where $\tau=\sigma^{\prime}(\alpha)$, $\rho=\sigma^{\prime}(\beta)$, $\rho >\tau$.

\underline{\bf Lemma 3.1.}

a). The equations ({\bf 3}.15) have complex roots with negative real part if and only if:
\begin{equation}
\varepsilon >0 \; \mbox{and} \; \varepsilon <\frac{\sqrt{\tau}}{\displaystyle 
nsin\frac{a_{r}}{2}},\; r=max\; J(\alpha)
\end{equation}

The equations ({\bf 3}.15) have real negative roots if and only if:
\begin{equation}
\varepsilon >0 \; \mbox{and} \; \varepsilon >\frac{\sqrt{\tau}}{\displaystyle
nsin\frac{a_{q}}{2}},\; q=min\; J(\alpha)
\end{equation}

b). The equations ({\bf 3}.16) have complex roots with negative real part if and only if:
\begin{equation}
\varepsilon >0 \; \mbox{and} \; \varepsilon <\frac{\sqrt{\rho}}{\displaystyle
nsin\frac{a_{r}}{2}},\; r=max\; J(\beta)
\end{equation}

The equations ({\bf 3}.16) have real negative roots if and only if:
\begin{equation}
\varepsilon >0 \; \mbox{and} \; \varepsilon >\frac{\sqrt{\rho}}{\displaystyle
nsin\frac{a_{q}}{2}},\; q=min\; J(\beta)
\end{equation}

c). The equations ({\bf 3}.17) have complex roots with negative real part if and only if:
\begin{equation}
\varepsilon >0 \;\frac{\sqrt{\rho+\tau-\sqrt{\Delta}}}{\displaystyle nsin\frac{a_{q}}
{2}}< \varepsilon <\frac{\sqrt{\rho+\tau+\sqrt{\Delta}}}{\displaystyle nsin\frac{a_{r}}{2}}
\end{equation}
$$q=min\; K(\alpha,\beta)),\; r=max\; K(\alpha,\beta))$$
where
$$\Delta =(\rho +\tau)^{2}-(\rho-\tau)^{2}cos^{2}\frac{a_{p}}{2}$$

\underline{\it Proof: } From the equations ({\bf 3}.15) and ({\bf 3}.16) it follows
directly a). and b).

The equations ({\bf 3}.17) may be written in the following form:
\begin{equation}
\lambda^{2}+4\varepsilon n^{2} sin^{2} \frac{a_{p}}{2}\lambda +2n^{2}(\rho +\tau)
sin^{2}\frac{a_{p}}{2}-in^{2}(\rho -\tau)sin\; a_{p}=0, p\in K(\alpha,\beta)
\end{equation}

The necessary and sufficient condition for these equations to admit complex roots
with negative real part is:
\begin{equation}
2n^{2}(\rho+\tau)sin^{2}\frac{a_{p}}{2} >\frac{16\varepsilon n^{4}sin^{4} \displaystyle
\frac{a_{p}}{2}}{8}+2\frac{n^{4}(\rho-\tau)^{2}sin^{2}a_{p}}{16\varepsilon^{2}n^{4}
\displaystyle sin^{4}\frac{a_{p}}{2}}
\end{equation}
and
$$4\varepsilon n^{2}sin^{2}\frac{a_{p}}{2}>0,\; 2n^{2}(\rho+\tau)sin^{2}\frac{a_{p}}{2}>0$$

\underline{\bf Theorem 3.5.} The 2-phase steady solution is stable if one of the
following conditions holds:

a).
$$\frac{\sqrt{\rho+\tau-\sqrt{\Delta}}}{\displaystyle nsin\frac{a_{q}}
{2}}< \varepsilon <\frac{\sqrt{\tau}}{n}, \; q=min\; K(\alpha,\beta)) $$

b).
$$\frac{\sqrt{\rho+\tau-\sqrt{\Delta}}}{\displaystyle nsin\frac{\pi}
{2n}}< \varepsilon <\frac{\sqrt{\rho+\tau+\sqrt{\Delta}}}{\displaystyle nsin\frac{a_{r}}{2}},
\; r=max\; K(\alpha,\beta)) $$

c).
$$max\left( \frac{\rho}{\displaystyle nsin\frac{a_{q_{2}}}{2}},\frac{\sqrt{\rho+\tau-
\sqrt{\Delta}}}{\displaystyle nsin\frac{a_{q_{3}}}{2}}\right) <\varepsilon <
min\left( \frac{\sqrt{\tau}}{\displaystyle nsin\frac{a_{r_{1}}}{2}},\frac{\sqrt{\rho+
\tau+\sqrt{\Delta}}}{\displaystyle nsin\frac{a_{r_{3}}}{2}}\right) $$
$$q_{2}=min\; J(\beta),\; q_{3}=min\; K(\alpha,\beta))$$
$$r_{1}=max\; J(\alpha),\; r_{3}=max\; K(\alpha,\beta))$$

d).
$$max\left( \frac{\sqrt{\tau}}{\displaystyle nsin\frac{a_{q_{1}}}{2}},\frac{\sqrt{\rho+\tau-
\sqrt{\Delta}}}{\displaystyle nsin\frac{a_{q_{3}}}{2}}\right) <\varepsilon <
min\left( \frac{\rho}{\displaystyle nsin\frac{a_{r_{2}}}{2}},\frac{\sqrt{\rho+
\tau+\sqrt{\Delta}}}{\displaystyle nsin\frac{a_{r_{3}}}{2}}\right) $$
$$q_{1}=min\; J(\alpha),\; q_{3}=min\; K(\alpha,\beta))$$
$$r_{2}=max\; J(\beta),\; r_{3}=max\; K(\alpha,\beta))$$

The proof is a direct consequence of the {\it Lemma 3.1.}

\underline{\bf Theorem 3.6.} If $\varepsilon =0$ the curves
$$u^{c}_{p}(t)=\bar{u}_{p}+exp\; (\lambda^{c}_{p}t) \sum_{k=1}^{n} \eta_{k}exp\;
(ia_{p}k), \; p\in J(\alpha)$$
$$u^{c}_{p}(t)=\bar{u}_{p}+exp\; (\lambda^{c}_{p}t) \sum_{k=1}^{n} \eta_{k}exp\;
(ia_{p}k), \; p\in J(\beta)$$
where $\lambda^{c}_{p},\; p\in J(\alpha)$ respectively $\lambda^{c}_{p},\; p\in J(\beta)$
are roots of the equations:
$$\lambda^{2}+4n^{2}sin^{2}\frac{a_{p}}{2}\tau =0,\; p\in J(\alpha) $$
$$\lambda^{2}+4n^{2}sin^{2}\frac{a_{p}}{2}\rho =0,\; p\in J(\beta) $$
are periodic with the period
$$T_{p}=\frac{\pi}{n\sqrt{\tau}} cosec \frac{p\pi}{2n}, \; p\in J(\alpha)$$
respectively
$$T_{p}=\frac{\pi}{n\sqrt{\rho}} cosec \frac{p\pi}{2n}, \; p\in J(\beta)$$

The proof follows from the above theorem.

\underline{\bf Remark.} For $\varepsilon<0$ the equations ({\bf 3}.15) and ({\bf 3}.16)
have positive roots and it follows that the 2-phase steady solution is unstable.

For $\varepsilon<0$ or $\varepsilon$ not satisfying one of the above conditions a).,
b)., c)., the roots of the equations ({\bf 4}.10), ({\bf 4}.11), ({\bf 4}.12)
may have positive or negative real part in the case they are complex, or they
are positive or negative real numbers. In this case the steady 2-phase solution
is unstable; it is hyperbolic saddle. The associated linear manifolds may be described
completely in the same manner we have done for the uniphase steady solutions.

\setcounter{equation}{0}
\begin{center}
\Large \bf 4. The discretization of the equation $u_{tt}=(\sigma (u_{x}))_{x}+
\varepsilon u_{xxt}$
\end{center}

\vspace*{.5cm}

Let $x_{k}=kh_{1}$, $k=1,...,n$ the division points of the interval $[0,1]$ with
$h_{1}=1/n$ and $t_{p}=ph_{2}$, $p=1,...,m$ the division points of the interval
$[0,1]$, $h_{2}=1/n$.

We$^{\prime}$ll approximate the derivates $\dot{u}_{k}(t)$, $\ddot{u}_{k}(t)$ by
$\displaystyle \frac{1}{h_{2}}[u_{k}^{p}-u_{k}^{p-1}]$ and $\displaystyle \frac{1}
{h_{2}^{2}}[u_{k}^{p+1}-2u_{k}^{p}+u_{k}^{p-1}]$, where $u_{k}^{p}=u(kh_{1},ph_{2})$.

We call the discrete system associated to the equation ({\bf 1}.1) the following
system:
$$\frac{1}{h_{2}^{2}} [u_{k}^{p+1}-2u_{k}^{p}+u_{k}^{p-1}]-\frac{1}{h_{1}} \sigma
\left( \frac{1}{h_{1}} (u_{k+1}^{p}-u_{k}^{p})-\frac{1}{h_{1}} \sigma (\frac{1}{h_{1}}
(u_{k}^{p}-u_{k-1}^{p})\right) =$$
\begin{equation}
=\frac{\varepsilon}{h_{1}^{2}h_{2}} [u_{k+1}^{p+1}-u_{k+1}^{p}-2(u_{k}^{p+1}-
u_{k}^{p})+u_{k-1}^{p+1}-u_{k-1}^{p}]
\end{equation}
$\hspace*{5cm} k=1,...,n-1,\; p=1,...,m-1$

The corresponding boundary conditions are:
\begin{equation}
u_{0}^{p}=0,\; u_{n}^{p}=P>0,\; p=1,...,m
\end{equation}

The system ({\bf 4}.1) represent the discrete Euler-Lagrange equations for the discrete
Lagrange function:
\begin{equation}
L(u_{k-1}^{p},u_{k}^{p},u_{k}^{p+1})=\frac{1}{2h_{2}^{2}} (u_{k}^{p+1}-u_{k}^{p})^{2}+
w(\frac{1}{h_{1}} (u_{k}^{p}-u_{k-1}^{p}))
\end{equation}
where $w:{\bf R}\to {\bf R}$, $w^{\prime}(\xi)=\sigma(\xi)$, and the dispersive term is:
$$\frac{\varepsilon}{h_{1}^{2}h_{2}} [u_{k+1}^{p+1}-u_{k+1}^{p}-2(u_{k}^{p+1}-
u_{k}^{p})+u_{k-1}^{p+1}-u_{k-1}^{p}]$$

For ({\bf 4}.1) it follows that the steady solutions $\bar{u}=\bar{u}_{k}$ satisfy:
\begin{equation}
\sigma (\frac{1}{h_{1}} (\bar{u}_{k+1}^{p}-\bar{u}_{k}^{p}))=\sigma (\frac{1}{h_{1}}
(\bar{u}_{k}^{p}-\bar{u}_{k-1}^{p})),\; k=1,...,n-1
\end{equation}

From ({\bf 4}.2) we obtain:
$$\sum_{k=1}^{n} \frac{1}{h_{1}} (\bar{u}_{k}^{p}-\bar{u}_{k-1}^{p})=\frac{1}{h_{1}}P $$

Thus a steady solution $\bar{u}$ satisfies the following:
\begin{equation}
\bar{u}_{k}^{p}-\bar{u}_{k-1}^{p}=P-\bar{u}_{n-1}^{p}, \; \sigma (\frac{1}{h_{1}}
(\bar{u}_{k}^{p}-\bar{u}_{k-1}^{p}))=C,\; k=1,...,n,\; C\in {\bf R}
\end{equation}

\underline{\bf Definition 4.1.} A steady solution $\bar{u}=(\bar{u}_{k}^{p})$ with
$\bar{u}_{k}^{p}=kh_{1}P$, $k=1,...,n-1$, $p=1,...,m$ is called uniphase steady solution.

A steady solution $\bar{u}=(\bar{u}_{k}^{p})$ with the property $\bar{u}_{k}^{p}-
\bar{u}_{k-1}^{p}=h_{1}f(k)$, $k=1,...,n-1$ is called multiphase steady solution.

Let $E=\{ \bar{u}=(\bar{u}_{k}^{p}) \} $ the set of multiphase steady solution. We
consider the following sets:
$$E^{+}=\left\{ \bar{u}=(\bar{u}_{k}^{p}):\sigma^{\prime} (\frac{1}{h_{1}}(\bar{u}_{k}^{p}-
\bar{u}_{k-1}^{p}))>0,\; k=1,...,n-1 \right\}$$
$$E^{-}=\left\{ \bar{u}=(\bar{u}_{k}^{p}):\sigma^{\prime} (\frac{1}{h_{1}}(\bar{u}_{k}^{p}-
\bar{u}_{k-1}^{p}))>0,\; k=1,...,n-1 \right\}$$

\underline{\bf Theorem 4.1.} The system ({\bf 4}.1) has a finite number of multiphase
steady solutions.

Let $\bar{u}=(\bar{u}_{k}^{p})$ a steady solution of ({\bf 4}.1) and $\tilde{u}_{k}^{p}=
\bar{u}_{k}^{p}+\displaystyle \sum_{l=1}^{n} \eta_{l}w_{l}^{p}(k)$, $k=1,...,n-1$
the components of a vector with $\eta_{l} \in (-a,a)$.

The linearized system associated to the system ({\bf 4}.1) in a neighborghood
of the steady solution $\bar{u}$ is:
$$h_{1}^{2}(w_{l}^{p+1}-2w_{l}^{p}+w_{l}^{p-1})-h_{2}^{2}h_{1}^{2}\left[ 
\sigma^{\prime}(\frac{1}{h_{1}} (\bar{u}_{k+1}^{p}-\bar{u}_{k}^{p})(w_{l+1}^{p}
-w_{l}^{p})-\right. $$
$$\left. \sigma^{\prime}(\frac{1}{h_{1}} (\bar{u}_{k}^{p}-\bar{u}_{k-1}^{p})
(w_{l}^{p}-w_{l-1}^{p})\right] =$$
\begin{equation}
=\varepsilon h_{2} [w_{l+1}^{p+1}+w_{l+1}^{p}-2(w_{l}^{p+1}-w_{l}^{p})+w_{l-1}^{p+1}-
w_{l-1}^{p}]
\end{equation}
$\hspace*{5cm} p=1,...,m,\; l=1,...,n-1$

A solution of ({\bf 4}.6) has the form:
\begin{equation}
w_{l}^{p}=\lambda^{p} exp\; (ia_{k}l),\; a_{k}=\frac{4\pi}{n},\; k=1,...,n,\; \lambda \in {\bf C}
\end{equation}

Replaing in ({\bf 4}.6) we obtain:

$$h_{1}^{2}(\lambda^{2}-2\lambda+1)exp\; (ia_{k})-h_{2}^{2}h_{1}^{2}(exp\; 2ia_{k}-
exp\; ia_{k})-$$
\begin{equation}
-s_{k}\lambda(exp\; ia_{k}-1)-\varepsilon h_{2}(\lambda^{2}-\lambda)
(exp\; ia_{k}-1)^{2}=0
\end{equation}

$\hspace*{5cm} k=1,...,n$ \\
where
$$s_{k}=\sigma^{\prime} \left( \frac{1}{h_{1}} (\bar{u}_{k}^{p}-\bar{u}_{k-1}^{p})\right) $$

From ({\bf 4}.8) it follows:
$$\lambda^{2}[h_{1}^{2}-\varepsilon h_{2}(exp\; ia_{k}-1)^{2}]-\lambda[2h_{1}^{2}
exp\; ia_{k}+h_{2}^{2}h_{1}(s_{k+1} exp\; ia_{k})(exp\; ia_{k}-1) -$$
\begin{equation}
-s_{k}(exp\; ia_{k}-1)-h_{2}\varepsilon (exp\; ia_{k}-1)^{2}]+h_{1}^{2}exp\; ia_{k}=0
\end{equation}

$\hspace*{5cm} k=1,...,n$

For a uniphase steady solution $\bar{u}=(\bar{u}_{k}^{p})$, $\bar{u}_{k}^{p}=
kh_{1}P$ and $\tau=\sigma^{\prime} (P)>0$, the equations ({\bf 4}.9) become:
\begin{equation}
\lambda^{2}(h_{1}^{2}-\varepsilon h_{2}\mu_{k})-\lambda(2h_{1}^{2}-h_{2}^{2}
h_{1}\tau \mu_{k}-\varepsilon h_{2}\mu_{k})+h_{1}^{2}=0
\end{equation}
$\hspace*{5cm} k=1,...,n$\\
where
$$\mu_{k}=\frac{(exp\; (ia_{k}-1))^{2}}{exp\; (ia_{k}-1)} =-4sin^{2}\frac{k\pi}{2n}$$

We suppose for the rest of this section that $h_{1},h_{2}\in (0,1)$ and $\varepsilon
\neq \displaystyle \frac{h_{1}^{2}}{h_{2}\mu_{k}}$.

\underline{\bf Theorem 4.1.} The uniphase steady solution for the system ({\bf 4}.1)
is assimptotically stable if and only if the following condition hold:
\begin{equation}
\tau >\frac{h_{1}}{h_{2}^{2}sin^{2}\frac{\pi}{2n}},\; \varepsilon \in \left(
-\frac{h_{1}^{2}}{2h_{2}sin^{2}\frac{\pi}{2n}}+\frac{h_{1}h_{2}}{2}\tau,\infty \right)
\end{equation}

\underline{\it Proof: } The equations ({\bf 4}.10) have the roots in modulus less
than 1 if and only if ({\bf 4}.11) hold.

A necessary and sufficient condition for the $k$-th equation in ({\bf 4}.10) to
have the modulus of its roots less than one is:
$$\frac{h_{1}^{2}}{h_{1}^{2}-\varepsilon h_{2}\mu_{k}}-1<0$$
\begin{equation}
-1-\frac{h_{1}^{2}-\varepsilon h_{2}\mu_{k}}{h_{2}^{2}}<\frac{|2h_{1}^{2}+h_{1}^{2}h_{1}
\tau \mu_{k}-\varepsilon h_{2}\mu_{k}|}{h_{1}^{2}}<1+\frac{h_{1}^{2}-\varepsilon
h_{2}\mu_{k}}{h_{1}^{2}}
\end{equation}

The inequatities ({\bf 4}.12) hold if and only if
\begin{equation}
\tau > -\frac{4h_{1}}{h_{2}^{2}\mu_{k}},\; \varepsilon \in \left( \frac{2h_{1}^{2}}
{\mu_{k}h_{2}}+\frac{h_{1}h_{2}}{2} \tau, \infty \right)
\end{equation}

But $\mu_{k}>\mu_{k+1}$ so it follows that ({\bf 4}.13) hold for any $k=1,...,n$
if and only if ({\bf 4}.11) hold.

\underline{\bf Theorem 4.2.} The uniphase steady solution $\bar{u}=(\bar{u}_{k}^{p})$
of the system ({\bf 4}.1) is unstable if one of the following conditions hold:

a).
$$\tau < -\frac{h_{1}}{h_{1}^{2}sin^{2}\frac{\pi}{2n}},\; \varepsilon \in \left( 
-\infty ,-\frac{h_{1}^{2}}{4h_{2}sin^{2}\frac{\pi}{2n}}+\frac{h_{1}h_{2}}{2} \tau \right) $$

b).
$$\varepsilon \in \left( -\infty ,-\frac{h_{1}^{2}}{2h_{2}sin^{2}\frac{\pi}{2n}}-h_{1}h_{2} \tau \right) $$

c).
$$\tau > \frac{h_{1}}{h_{1}^{2}sin^{2}\frac{\pi}{2n}},\; \varepsilon \in \left( 
-\frac{h_{1}^{2}}{4h_{2}sin^{2}\frac{\pi}{2n}},-\frac{h_{1}^{2}}{4h_{2}sin^{2}
\frac{\pi}{2n}}+h_{1}h_{2}\tau \right) $$

\underline{\it Proof :} a). The equations ({\bf 4}.10) have the roots in modulus
greater than one only if a). hold. A necessary and sufficient condition for the
$k$-th equation from ({\bf 4}.9) to have roots in modulus greater than one is:
$$\frac{h_{1}^{2}-\varepsilon h_{2}\mu_{k}}{h_{1}^{2}}-1 <0$$

\begin{equation}
-1-\frac{h_{1}^{2}-\varepsilon h_{2}\mu_{k}}{h_{1}^{2}}<\frac{|2h_{1}^{2}+h_{2}^{2}h_{1}
\tau \mu_{k}-\varepsilon h_{2}\mu_{k}|}{h_{1}^{2}}<\frac{h_{1}^{2}-\varepsilon
h_{2}\mu_{k}}{h_{1}^{2}}+1
\end{equation}

From ({\bf 4}.14) we obtain
\begin{equation}
\tau < -\frac{4h_{1}}{h_{2}^{2}\mu_{k}},\; \varepsilon \in \left( -\infty ,\frac{2h_{1}^{2}}
{h_{2}\mu_{k}}+\frac{h_{1}h_{2}}{2} \tau \right)
\end{equation}

b)., c). can be prooved in an analogous manner.

\underline{\bf Theorem 4.3.} The uniphase steady solution $\bar{u}=(\bar{u}_{k}^{p})$
for the system ({\bf 4}.1) with $\varepsilon =0$ is unstable if
$$\tau <\frac{h_{1}}{h_{2}^{2}sin^{2}\frac{\pi}{2n}} $$

\underline{\it Proof :} For $\varepsilon=0$, the equation ({\bf 4}.10) become:
\begin{equation}
\lambda^{2}-\lambda (2+\frac{h_{2}^{2}}{h_{1}}\tau \mu_{k})+1=0 ,\; k=1,...,n
\end{equation}

From ({\bf 4}.16) we have $\lambda_{1}\lambda_{2}=1$ and $\displaystyle 
\lambda_{1}+\lambda_{2}=2+\frac{h_{2}^{2}}{h_{1}}\tau \mu_{k}$. A necessary and
sufficient condition for $k$-th equation from ({\bf 4}.16) to have one root greater
than one and the other less than one is:
$$\tau <\frac{h_{1}}{h_{2}^{2}sin^{2}\frac{\pi}{2n}}$$

\underline{\bf Proposition 4.1.} Let $u=\bar{u}_{k}^{p}$ a uniphase steady solution
for the system ({\bf 4}.1) which satisfies one of the condition b). or c). from the
theorem 4.2. and $\lambda_{k}^{+}$, $\lambda_{k}^{-}$ the roots of the $k$-th
equation ({\bf 4}.9) with $|\lambda_{k}^{+}|>1$ and $|\lambda_{k}^{-}|<1$.

The linear stable manifold is given by:
$$u_{k}^{sp}=kh_{1}P+(\lambda^{-}_{k})^{p} \sum_{l=1}^{n} \eta_{l}exp\; (ia_{k}l)$$
$\hspace*{5cm} p=1,...,m,\; k=1,...,n-1$

The linear unstable manifold is given by:
$$u_{k}^{up}=kh_{1}P+(\lambda^{+}_{k})^{p} \sum_{l=1}^{n} \eta_{l}exp\; (ia_{k}l)$$
$\hspace*{5cm} p=1,...,m,\; k=1,...,n-1$

In these conditions it follows that:
$$\lim_{p\to \infty} u_{k}^{sp}=kh_{1}P \; \mbox{and} \; \lim_{p\to -\infty} u_{k}^{up}=kh_{1}P$$

\begin{center}
\bf References
\end{center}

[1]. S. Bir\u{a}ua\c{s}, S. Balint, A. M. Balint - {\it On an infinite system
of differential equations related to the hydrodynamical shock problems},
Bulletins for Appl. Math. (BAM), 1260/96 (LXXX), pp. 247-259\\

[2]. J. M. Greenberg - {\it Continuum limite of discret gases}, Arch. Rat. Mech.
Anal. 105 (1989), pp. 367-376\\

[3]. P. D. Lax - {\it On dispersive difference schemes}, Phisica D, 1986, pp. 250-254\\

[4]. J. von Neumann - {\it Proposal and analysis of a numerical method for the
treatment of hydrodynamical shock problems}, Collected Work, vol. VI, pp. 361-379,
Pergamon, Oxford, 1963

\end{document}